\documentstyle[12pt]{amsart}
\begin{document}
\newcommand{\ep}{exp(\frac{-1}{\pi} \int_{\Omega} \frac{dA(\zeta)}
                   {(\zeta - z)({\overline{\zeta}}-{\overline{w}})} )}
\newcommand{\Resw}{({T^{\ast}}-{\overline{w}}{)^{-1}}\xi}
\newcommand{\Resz}{({T^{\ast}}-{\overline{z}}{)^{-1}}\xi}
\newcommand{\su}{ sup \{ \frac{{l_a}(p)}{ \| p {\|_{1,K}} } ; 
p \in {{\cal P}_N} \setminus \{0\} \} }
\newcommand{\inte}{ \frac{d{\overline{\zeta}}\wedge d\zeta}
                {(\zeta -z)(\overline{\zeta} -\overline{w})}  }
\newcommand{\binte}{ \frac{log(\overline{\zeta}-\overline{w})d\zeta}
                 {\zeta -z} }
\newcommand{\variabile}{ {z_1},\overline{z_2};{w_1},\overline{w_2}}
\newcommand{\presu}{(T-u{)^{-1}}}
\newcommand{\presv}{(T-v{)^{-1}}}
\newcommand{\resu}{({T^\ast}-\overline{u}{)^{-1}}}
\newcommand{\presza}{(T-{z_1}{)^{-1}}}
\newcommand{\reszb}{({T^\ast}-\overline{z_2}{)^{-1}}}
\newcommand{\preswb}{(T-{w_2}{)^{-1}}}
\newcommand{\reswa}{({T^\ast}-\overline{w_1}{)^{-1}}}
\newcommand{\PRESwa}{(T-{w_1}{)^{-1}}}
\newcommand{\RESwb}{({T^{\ast}}-\overline{w_2}{)^{-1}}}

\def\paragraph{\subsection*}

\title{Extremal solutions of the two-dimensional $L$-problem of
moments, II}
\author{Mihai Putinar}
\address{Department of Mathematics, University of California, 
Riverside, CA 92521}
\email{mputinar@@ucrmath.ucr.edu}
\thanks{This paper was completed while the author was 
on leave from University of California at Riverside, visiting
the Mathematical Sciences Research Institute in Berkeley; paper 
supported in part by NSF grant DMS 9500954 and
research at MSRI partially supported by NSF grant DMS 9022140.}

\begin{abstract}
All extremal solutions of the truncated $L$-problem of moments in two real
variables , with support contained in a given compact set,
are described as characteristic functions of semi-algebraic sets given by
a single polynomial inequality. An exponential kernel, arising as the
determinantal function of a naturally associated hyponormal operator with
rank-one self-commutator, provides a natural defining function for these
semi-algebraic sets. We find an intrinsic characterization of this 
kernel and we describe a series of analytic continuation properties of it
which are closely related to the behaviour of the Schwarz reflection 
function in portions of
the boundary of the extremal supporting set.
\end{abstract}

\maketitle

\section{Introduction}  In a previous paper, [P1], a special class of extremal solutions
of the $L$-problem of moments in two variables was related via hyponormal
operators to quadrature domains in the plane and to some rational functions 
related to them. It is the aim of the present paper to apply the same techniques
to all extremal solutions of the $L$-problem and to begin a study of the analytic
objects arising from this investigation. Although this programme is more
general and and the results below are less precise, a series of facts from
operator theory and function theory converge to 
an interesting picture of all extremal solutions of the $L$-problem of moments.

To explain the above sentences we have to  be more specific. Let $K$ be a 
compact subset of the complex plane, let
$L$ be a fixed positive constant and let $N$ be a fixed positive integer. 
We are interested in classifying and characterizing in intrinsinc terms the
moments $${a_{kl}}= \int_{K} \phi(x,y) {x^n}{y^m} dA,\hspace{.2in} 0 \leq m \leq m+n \leq N,
$$ 
of a measurable function $\phi$ on $K$ which satisfies $0 \leq \phi \leq L$, $dA$-a.e. .
Above, and throughout the paper $dA$ stands for the planar Lebesque measure.
Note that, because we are working with the two-dimensional Lebesgue measure, the
alterations of null measure of the compact set $K$ (such as adding or removing 
continuous exterior lines or
internal slits) are 
not significant for our discussion.
As a matter of fact, later in the paper we will restict ourselves to the case when $K$ is
the closure of a domain with real algebraic boundary.

Let  $\Sigma$ denote the collection of all vectors $a=({a_{kl}}{)_{k+l \leq N}} \in
{{\bf R}^d},(d= \frac{(N+1)(N+2)}{2} )$ which arise as the moments of a function
$\phi$ as before. It is clear that $\Sigma$ is a compact convex subset of
${{\bf R}^d}$. Following M.G.Krein and his classical by now convexity theory (cf.
[K],[KN])
we infer that a point $a$ is extremal in the set $\Sigma$ if and only if it 
corresponds to the moments of a function $\psi$ of the form $$
\psi = L {\chi _ \Omega},\hspace{.2in} \Omega= \{ (x,y) \in K ; p(x,y) > 0 \}, $$
where $p$ is an arbitrary polynomial of degree $N$ (with real coefficients)
and ${\chi _D}$ is the characteristic function of the set $D$. Moreover, we will
see below that only in this case the above moment problem has a unique 
solution. The role of the bound $L$ in some related extremality problem will also
appear in the sequel. Thus, following the classical one-variable theory, we may ask how the
extremal solution $\phi$ is encoded in its moments of degree less or equal than
$N$. Although we are far from having a satisfactory answer to this basic question
in its full generality, some particular cases are worth being discussed in more detail.

For instance, in the special case when the set $\Omega$ above (in the formula of the extremal
solution $\phi$) is a quadrature domain, the following
exponential kernel was the key to the preceding determination problem (and actually 
to much more):
$$
{E_{\Omega}}(z,{\overline{w}})=\ep, \hspace{.2in} (z,w \in {\bf C} \setminus {\overline{\Omega}}).
$$
(The reader will notice that we have passed tacitly to complex coordinates; this
transformation obviously does not change the moment problem). 

The importance of the
exponential kernel is two fold: first it is analytic in $z,{\overline{w}}$ and the
moments of the function ${\chi _{\Omega}}$ can be deduced by simple algebraic operations
from its Taylor expansion at infinity, and second it admits a canonical factorization as
$$ {E_{\Omega}}(z,{\overline{w}})=1-\langle \Resw , \Resz \rangle,\hspace{.2in} 
(z,w \in {\bf C} \setminus {\overline{\Omega}}).$$
where $T$ is the unique hyponormal operator with one-dimensional commutator
$$( [{T^{\ast}},T] =\xi \otimes \xi)$$ having ${\chi _{\Omega}}$ as principal function.
See [P1],[P2] for details and bibliographical references.

The main result of [P1] asserts that $\Omega$ is a domain of quadrature
(in the sense of [Gu],[Sa],[Sh]) if and only if there is a polynomial $P(z)$
such that the function $P(z){\overline{P(w)}}{E_{\Omega}}(z, {\overline{w}})$
is polynomial at infinity,or, if and only if the ${T^{\ast}}$-invariant subspace
generated by the vector $\xi$ is finite dimensional. In this case a simple and rather
constructive dictionary relates the above three objects, see [P2]. For
what follows it is important to remark that, in the case of a quadrature domain,
the associated exponential kernel extends meromorphically in each variable 
from ${\bf C} \setminus {\overline{\Omega}}$ to the whole plane. As explained
in [P2], in the polar parts of the Laurent series of ${E_{\Omega}}$ at its
(finitely many) poles we can read the complete sequence of moments of 
${\chi _{\Omega}}$,
the operator $T$ and in particular a defining equation for the domain $\Omega$.

Let us turn now to an arbitrary domain $\Omega \subset K$ which carries an extremal
solution of the truncated $L$-problem of moments supported by $K$. A first part of
the present paper is devoted to the analytic continuation properties of the
associated exponential kernel ${E_{\Omega}}(z,{\overline{w}})$. First we will see
that ${E_{\Omega}}$ extends analytically in both variables $z, {\overline{w}}$
across a point $\lambda \in \partial \Omega$ only if the boundary $\partial \Omega$
is real analytic in a neighbourhood of $\lambda$. Then we relate the local
analytic continuation of the exponential kernel to the Schwarz function of
the boundary of $\Omega$. (See [D],[Sh] for details about the Schwarz function).

Passing from local to global analytic continuation, we prove that the kernel
${E_{\Omega}}$ extends analytically inside $\Omega$ as far as the Schwarz function
(of a portion of the boundary ) extends along continuous paths. 

Suppose now that the supporting compact set $K$ is semi-algebraic, that is it is defined
by a simultaneous system of polynomial inequalities. Such examples are the disk,
the square, etc. Thus, in view of the aforementioned characterization of 
extremal solutions $\phi$ of the $L$-problem of moments, the supporting 
domain $\Omega$ of $\phi$ is still semi-algebaric (with one more defining inequality).
In that case it is well known that the Schwarz function of the irreducible components
of the boundary of $\Omega$ is an algebraic function. Thus the kernel ${E_{\Omega}}$
extends in each variable, along paths, from a connected component $C$ of $\partial \Omega$ 
to the whole
domain $\Omega$, minus the ramification points of the respective Schwarz 
functions. As a consequence we can produce quadrature formulas for the domain
$\Omega$, supported as anaytic functionals on a system of curves (more
precisely the cuttings which specify a determination of the multivalued
Schwarz function).

We would like to mention that, as in the preceding papers, the operator theory
is necessary in several points of the subsequent proofs, although the statements are
purely function theoretic.

The paper is organized as follows. Section 2 recalls the essential results of
M.G.Krein which characterize the extremal solutions of the truncated $L$-problem of
moments. A great deal of effort was put in this area by statisticians, who have
followed A.A.Markov and P.L. Chebyshev in deriving sharp estimates for
the distribution of collections of random variables with a set of moments 
prescribed, see [G],[KS]. Section 3 is devoted to local and global analytic 
continuation properties of the exponential kernel of a planar domain. In Section
4 we specialize the results to semi-algebraic domains and we obtain a general
quadrature formula supported on thin sets. Section 5 is independent of the main
body of the paper; we describe there a set of positivity conditions which characterize
the exponential kernel. In this section we meet old ideas and techniques 
(based on the magic of resolvents of linear operators) due
to de Branges and his followers. In particular, we interpret in this final section
the $L$-problem of moments as an interpolation problem of the Carath\'{e}odory-Fej\'{e}r
type, for a class of analytic functions defined in a polydisk of
dimension four.\bigskip

{\it Acknowledgements.} We would like to thank the Mathematical Sciences 
Research Institue in Berkeley for creating a unique scientific environment
which made possible this work. 
We are indebted to Bj\"{o}rn Gustafsson and Harold S.Shapiro for several
discussions on related themes and for their interest in this subject. We
would like to thank James Rovnyak for his comments on the
structure of determinantal functions and the bibliogarphical reference
[PR].

\section{Convexity results} 

This first part of the paper contains a survey of a series of known results 
derived from the work of M.G.Krein and his succesors. They are intended to serve as a 
motivation for the next sections. For that reason this part is independent of the rest
of the paper, and it is presented in a slightly more general setting (${{\bf R}^n}$
rather than ${{\bf R}^2}$ and an arbitrary support compact $K$ for the moment problem).
The basic monographs we refer to for details are [KS] and [KN].

The variable in ${{\bf R}^n}$ is denoted by $x=({x_1},\ldots ,{x_n})$; $dx$ means
the volume measure in ${{\bf R}^n}$. For a multi-index $\alpha =({\alpha _1},\ldots ,
{\alpha _n}) \in {{\bf N}^n}$ we denote $|\alpha |={\alpha _1}+\ldots +
{\alpha _n}$ and put as usual ${x^{\alpha}}={{x_1}^{\alpha _1}}\ldots
{{x_n}^{\alpha _n}}.$

Let $K$ be a compact subset of ${{\bf R}^n}$ ; in order to avoid some minor complications
we assume that $int(K) \neq \emptyset$. Fix a positive integer $N$ and a positive
constant $L$. The {\it truncated $L$-problem of moments} supported by the set
$K$ consists in finding necesary and sufficient conditions for a sequence
$a=({a_{\alpha}}{)_{|\alpha | \leq N}}$ of real numbers to be represented as:
\begin{equation}
{a_{\alpha}}= \int_{K} {x^{\alpha}} \phi(x) dx, \hspace{.2in}(|\alpha | \leq N,
\phi \in {L^1}(K,dx), 0 \leq \phi \leq L).
\end{equation}
Moreover, it is traditionally of interest to classify all solutions of this problem
and to characterize the uniqueness cases, see [KN] Chapter VII.

For a first part of this section we can adopt the normalization $L=1$.

Let us denote, for $L=1$, the set of all possible moment sequences as follows:
\begin{equation}
\Sigma = \{ a(\phi)=({a_{\alpha}}) ; {a_{\alpha}}=\int_{K} {x^{\alpha}}\phi (x) dx,
           | \alpha| \leq N, \phi \in {L^1}(K,dx), 0 \leq \phi \leq 1 \}.
           \end{equation}
Let ${\bf R}[x]$ be the space of polynomials in the variables $x$, with real coefficients.
We put ${{\cal P}_N}={{\cal P}_N}({{\bf R}^n})= \{ p \in {\bf R}[x], deg(p) \leq N \}$.
With these data fixed, we denote $d=dim({{\cal P}_N})$ and parametrize the coordinates
in the space ${{\cal P}_N} \subset {{\bf R}^d}$ as follows $y=({y_{\alpha}}{) 
_{|\alpha | \leq N }}$.

It is clear that $\Sigma$ is a compact convex subset of ${{\bf R}^d}$. Let ${a^0}$
be a boundary point of $\Sigma$ and let $f(y)=\langle c, y \rangle +d$ be a
supporting affine functional to $\Sigma$, at the point ${a^0}$. Let us represent
the point ${a^0}$ as the moment sequence of the function ${\phi _0}$:
${a^0}=a({\phi _0})$. Then we have the following relations: $$
\langle c,a \rangle +d \leq 0, \hspace{.2in}(a \in \Sigma),$$
and $$
\langle c, {a^0} \rangle +d =0.$$
By substarcting them and representing $a$ as $a=a(\phi)$ we obtain:
\begin{equation}
\int_{K} {p^0}(x)(\phi (x) -{\phi _0}(x)) dx \leq 0, \hspace{.2in}
(\phi \in {L^1}(K), 0 \leq \phi \leq 1), \end{equation}
where ${p^0}(x) = \sum_{|\alpha |\leq N} {c_\alpha}{x^{\alpha}}$. But
relation (3) is possible for all functions $\phi$ as above if and only if
${p^0}(x) >0$ implies ${\phi _0}(x) =1$ and ${p^0}(x)<0$ implies ${\phi _0}(x)
=0$. Since the set of zeroes of a non-trivial polynomial has null volume measure,
these latter implications determine an unique element ${\phi _0} \in {L^1}$. In other 
terms we have proved that ${\phi _0}= {\chi_{\{ {p^0}>0\} }}$ a.e. ,where
${\chi_S}$ is the charactersistic function of the set $S$.

In fact the above argument can easily be reversed, and we can state the following
conclusion.

\paragraph{Theorem 2.1.} {\it A point $({a_\alpha}{)_{|\alpha| \leq N}}$ belongs to the
boundary of the set $\Sigma$ of all moments if and only if there is a non-trivial
polynomial $p(x)$ of degree less or equal than $N$, and with the property:$$
{a_\alpha} = \int_{K \cap \{ p>0 \} } {x^\alpha} dx, \hspace{.2in}(|\alpha | \leq N). $$}\bigskip

Above we have denoted in short by $\{ p>0 \}$ the set of those points $x$
which satisfy $p(x)>0$. Since we have assumed the compact set $K$ to possess
interior points, it is immediate to remark that $a(\phi) \in int(\Sigma)$ for all
functions $\phi$ satisfying $0 < \phi <1$.

\paragraph{Theorem 2.2.(M.G.Krein)} {\it A necessary and sufficient condition
for the truncated 1-problem of moments (1) to be solvable is that, for every
polynomial $p(x)=\sum_{|\alpha | \leq N} {c_\alpha}{x^\alpha},$
we have:$$
\sum_{|\alpha | \leq N} {a_\alpha}{c_\alpha} \leq \int_{K} max(p(x),0) dx.$$}\bigskip

{\it Proof.}  Let us put $a=({a_\alpha}),c=({c_\alpha}), {p_{+}}=
max(p,0)$. The necessity follows from
the observation:$$
\langle a, c \rangle = \int_{K} p(x) dx \leq \int_{K} {p_{+}}(x) dx.$$

For proving the sufficiency we will show that the vector $a$ of virtual moments
and the point $0 \in \Sigma$ cannot be separated by a supporting hyperplane to
the set $\Sigma$. Exactly as before, let :$$
\langle c, y \rangle +d \leq 0, \hspace{.2in} (y \in \Sigma),$$
$$ \langle c, {y^0} \rangle +d =0, $$
be a supporting hyperplane to $\Sigma$ at the boundary point 
$${y^0}=a({\phi_0}), {\phi_0} =
{\chi _{{p^0} >0}}, {p^0}(x)= \sum {c_\alpha}{x^\alpha}.$$ In particular
$$ d= \langle c, {y^0} \rangle = -\int_{K} {{p^0}_{+}}(x) dx \leq 0.$$
Therefore, $ \langle c, 0 \rangle +d \leq 0$ and the proof is complete.\bigskip

By following the same lines one can obtain generalized Chebysev inequalities.
Namely, for a continuous function $\psi$ on $K$, not belonging to the space of
polynomials ${{\cal P}_N}$ and a point $a \in int(\Sigma)$ we want to find the
extreme values of $\int_{K} \psi \phi dx$ over all measurable functions 
$\phi, 0 \leq  \phi \leq 1,$  having the finite sequence of moments $a=a(\phi)$
prescribed. Let us denote $$
\Phi (a) =\{ \phi \in {L^1}(K); 0 \leq \phi \leq 1, a(\phi)=a \}.$$
The main result in this area, with important applications to mathematical statistics,
can be stated as follows: \bigskip

{\it There are polynomials ${\overline{p}},{\underline{p}} \in {{\cal P}_N}$
with the property that ${\chi_{\{\psi >{\overline{p}} \}}}, {\chi_{\{
{\underline{p}}>\psi \}}} \in \Phi (a)$
and $$
sup \{ \langle p, a \rangle ; p \in {{\cal P}_N}, p \leq \psi \} =
\int_{K \cap \{ {\underline{p}} >\psi \}} \psi dx = min_{\phi \in \Phi (a)} 
\int_K \psi \phi dx  \leq $$                                
$$ max_{\phi \in \Phi (a)} \int_{K} \psi \phi dx = \int_{K \cap \{ \psi > {\overline{p}} \} }
\psi dx = inf \{ \langle p,a \rangle ; p \in {{\cal P}_N}, p \geq \psi \}.$$ }   \bigskip

Above we identify the polynomial $p$ with the sequence of its coefficients, also
denoted by $p \in {{\bf R}^d}$.

For a proof and relevant comments about the above facts, see [KS] Sections
VIII.8 and XII.2. For applications of generalized Chebyshev inequalities we also
refer to [G].

For the aims of the present paper it is enough to retain from these inequalities
the fact that they are attained on some extremal solutions $\phi \in \Phi (a)$ always
given by the characteristic function of a sublevel set of a polynomial (plus possibly
a multiple of the new function $\psi$).

Now we slightly change the point of view. Since we assume the supporting compact set
given, after a translation and homothety the problem (1) is equivalent to:
$$ 2{a_\alpha} -L\int_{K} {x^\alpha} dx = \int_{K} {x^\alpha} (2\phi (x)-L)dx, \hspace{.2in}(|\alpha | \leq N),$$
where the new unknown function $2 \phi -L$ satisfies:$ -L \leq 2 \phi -L \leq L$.
Modulo this transformation we consider henceforth the following moment problem:
\begin{equation}
{a_\alpha} = \int_{K} {x^\alpha} \phi(x) dx,\hspace{.2in}(|\alpha |\leq N, \phi \in {L^1}(K), -L \leq \phi \leq L).
\end{equation}
Let us denote $a=({a_\alpha}{)_{|\alpha | \leq N}}$ to be the sequence of virtual moments.

Because we have assumed $int(K) \neq \emptyset$, the monomials $({x^\alpha}
{)_{|\alpha| \leq N}}$
are linearly independent regarded as functions of $x \in K$. Thus, for $L$ large
enough the problem (4) always admits a solution $\psi$.

Let us cosider the embedding ${{\cal P}_N}({{\bf R}^n}) \subset {{L^1}_{\bf R}}(K,dx)$ and
consider the linear functional $$
{l_a}:{{\cal P}_N} \longrightarrow {\bf R},\hspace{.2in} {l_a}({x^\alpha})={a_\alpha},
(|\alpha| \leq N).$$
In virtue of Riesz Theorem, any continuous extension of ${l_a}$ to 
${{L^1}_{\bf R}}(K)$ is represented by a function $\phi \in {{L^\infty}_{\bf R}}(K)$.
Hence $\phi$ is a solution of problem (4) and we have:
\begin{equation} 
{l_a}(p) \leq \| p{\|_{1,K}} \| \phi {\|_{\infty, K}},\hspace{.2in}(p \in {{\cal P}_N}).
\end{equation}

Then it remains to remark that the converse implication is given by Hahn-Banach Theorem. Moreover, 
the familiar analysis of the equality case in (5) is also relevant for us. Summing up, we can state the
next theorem.

\paragraph{Theorem 2.3.} {\it a). Problem (4) admits a solution if and only if
$$ L \geq \su .$$

b). If $L=\su$, then the solution is unique and it coincides with the function
$L sgn({p_0})$, where ${p_0}$ is a polynomial of degree less or equal than $N$.}\bigskip

Below we only sketch the proof of Theorem 2.3. The reader can find more details in
[KN] \S IX.1-2.

Let ${p_0} \in {{\cal P}_N} \setminus \{0\}$ be a polynomial with the property that
that $$
{L_0} = \su =\frac{ {l_a} ({p_0})} {\| {p_0} {\| _{1,K}} },$$
and let ${\psi _0} \in { {L^\infty}_{\bf R} }(K)$ be an extremal element which satisfies
${\psi _0}(p)={l_a}(p), p \in {{\cal P}_N}$, $\| {\psi_0} {\| _{\infty, K}}={L_0}$ and
$$ \int_{K} {p_0}{\phi_0} dx = \| {p_0} {\| _{1,K}} \| {\phi_0} {\| _{\infty ,K}}.$$

Since $\int_{K} ({p_0}{\phi_0}-|{p_0}| {L_0})dx=0$ and ${p_0}{\phi_0} \leq |{p_0}|{L_0}$
we obtain $$
{\phi_0}={L_0} sgn({p_0}), a.e. .$$
Notice again that the zero set of a polynomial has null measure and therefore
the element ${\phi_0}$ is well defined almost everywhere.

Moreover, if ${q_0}$ is anoter non-zero element of ${{\cal P}_N}$ satisfying
the extremal condition $\frac{ {l_a}({q_0})} { \| {q_0} {\|_{1,K}} }={L_0}$, then
$$
\int_{K} \frac{q_0}{ \| {q_0} {\|_{1,K}} } {\phi _0} dx =
\int_{K} \frac{p_0} { \| {p_0} {\| _ {1,K}}} {\phi_0} dx = {L_0}.$$
Whence ${\phi _0} ={L_0} sgn( {q_0}).$ 

Thus the extremal solution ${\phi_0}$ of the moment problem (4) is indeed unique.

>From the previous argument it is also clear that the problem (4) has a convex 
continuum of solutions whenever $L$ is greater than the critical value ${L_0}$.

To finish our brief presentation of these convexity methods, let us reverse
the problem (4) and present the above results in the following form.  

\paragraph{Corollary 2.4.} {\it The function 
$\phi \in {{L^\infty}_{\bf R}}(K)$ is uniquely determined in the ball
$\{ \psi  \in {{L^\infty}_{\bf R}}(K); \| \psi {\| _{\infty, K}} \leq \|
\phi {\|_{\infty, K}} \}$
by its finite sequence of moments $a=a(\phi)$
if and only if $$
\phi = \| \phi {\|_{\infty, K}}sgn(p),$$
where $p \in {{\cal P}_N} \setminus \{0\}$.

Conversely, for any non-zero polynomial $p \in {{\cal P}_N}$ the moemnts $$
{a_\alpha} = \int_{K \cap \{p>0\}} {x^\alpha} dx -
\int_{K \cap \{ p<0 \} } {x^\alpha}dx, \hspace{.2in}(|\alpha| \leq N),$$
determine the function $sgn(p)$ in the unit ball of ${{L^\infty}_{\bf R}}$.}\bigskip

For the latter statement, it suffices to remark that, for any function
$\phi \in {{L^{\infty}}_{{\bf R}}}(K)$ satisfying $a(\phi)=a(sgn(p))$ we have:
$$\| \phi {\|_{\infty ,K}} \| p {\|_{1,K}} \geq 
|\int_{K} \phi p dx | = | \int_{K} sgn(p) p dx| = \int_{K} sgn(p) p dx =
\| p{\|_{1,K}}.$$

If we assume in addition that $\| \phi {\|_{\infty, K}} =1$, then we obtain
the stated uniqueness result $\phi =sgn(p)$, a.e. .

To draw a conclusion of this section, we have seen from two different perspectives that
the extremal solutions of the moment problem (1) are parametrized by the 
sub-level sets $\{ p>0 \}$ of any non-zero polynomial $p$ of degree less or equal than
$N$. 

\section{Analytic continuation of the exponential kernel}

>From now on we restrict our study to the $L$-problem of moments in two real
variables. As explained in the previous paper [P1] in this dimension the theory
of hyponormal operators meets favorably the $L$-problem. Without entering into
all technical details exposed in [P1], [P2] we fix below some conventions 
and notation. An introduction to the theory of hyponormal operators is given in [MP].

Let $\phi$ be a measurable function with compact support in the complex plane
which satsifies $0 \leq \phi \leq 1$, a.e. .Let $T$ be the unique (up to unitary 
equivalence) hyponormal operator with rank-one self-commutator $([ {T^\ast},T]=
\xi \otimes \xi)$ with the principal function equal to $\phi$. 

The two objects above are related by the following formula:$$
\langle \Resw , \Resz \rangle = 1- exp( \frac{-1}{\pi} \int_{\bf C}
\frac{\phi(\zeta)dA(\zeta)}
{ (\zeta -z)({\overline{\zeta}}-{\overline{w}})} ),$$
which is valid for all points $z,w$ in the reslovent set of the operator $T$.
Actually a separately continuous extesion of the above formula to the whole 
${{\bf C}^2}$ holds, see [P1]. The importance of this formula lies in the
fact that it relates, after taking the Taylor expansions at infinity, the moments
of the function $\phi$ to the moments of the operator $T$.

We are interested in the sequel in the moments of the extremal solutions discussed
in the preceding section. Therefore, assuming the supporting compact set $K$ (introduced in \S 2)
nice, we will study the above relationship only for characteristic functions 
$\phi={\chi _{\Omega}}$ , where $\Omega$ is a bounded domain of the complex plane,
satisfying $\Omega = int({\overline{\Omega}})$. To simplify notation we put:
$$ {E_{\Omega}}(z,{\overline{w}})= \ep, \hspace{.2in}(z,w \in {\bf C} \setminus {\overline{\Omega}}).$$

This exponential kernel is analytic in $z$ and antianalytic in $w$. To simplify 
the terminology an analytic-antianalytic function of several variables
denoted as $$f(z,u,\ldots;{\overline{w}},{\overline{t}},\ldots),$$ will
be called analytic in $z,u,\ldots; {\overline{w}},{\overline{t}},\ldots$ .

To each domain $\Omega$ as above we associate the unique irreducible hyponormal operator
$T$ satisfying $[{T^\ast},T]=\xi \otimes \xi$ and having the principal function
equal to ${\chi _{\Omega}}$. We simply call $T$ {\it the hyponormal operator corresponding
to $\Omega$}. Recall that the spectrum of $T$ is the closure of $\Omega$, the
essential spectrum is the boundary of $\Omega$, and so on. See [MP] Chapter XI
for more details.

The aim of the present section is to investigate the analytic continuation
properties of the kernel ${E_{\Omega}}$ across portions of the boundary of
$\Omega$. This study is motivated by the earlier results obtained in the case
of quadrature domains, cf. [P1],[P2].

\paragraph{Proposition 3.1.} {\it Let $\lambda \in \partial \Omega$ and assume that there exists
an open neighbourhood $U$ of $\lambda$ in ${\bf C}$ with the property that
the function ${E_{\Omega}}(z,{\overline{w}})$ extends analytically from
$(U \setminus {\overline{\Omega}}) \times (U \setminus {\overline{\Omega}})$
to $U \times U$.

Then there is an open neighbourhood $V$ of $\lambda$ in $U$ such that
$V \cap \partial \Omega$ is a real analytic set, on which the extension
of the function ${E_{\Omega}}(z,{\overline{z}})$ vanishes.}\bigskip

{\it Proof.} Let ${E_{\Omega}}(z,{\overline{w}})$ be the canonical extension of the
exponential kernel to ${\bf C} \times {\bf C}$. We recall that, denoting by $\Resz$
the unique solution $\eta$ of minimum norm of the equation:$$
({T^\ast}-{\overline{z}})\eta =\xi, $$
we obtain a weakly continuous function defined everywhere on ${\bf C}$. 
Then the identity $$ \ep =1 -\langle \Resw , \Resz \rangle $$ holds everywhere, by
a result due to K.Clancey, see [P1].We also
know that $\| \Resz \| =1$ for all points $z \in \Omega$, cf. [MP] Theorem XI.4.1.

Let $F(z,{\overline{w}})$ be the analytic extension of ${E_{\Omega}}(z,{\overline{w}})$
to $U \times U$, whose existence was assumed in the statement.

Let us consider the closed linear span $K= \bigvee_{z \in U \setminus {\overline{\Omega}}}
\Resz$. Since the statement of Lemma 3.1 is local, we can
assume that $U \setminus {\overline{\Omega}}$ intersects the unbounded component 
of ${\bf C} \setminus {\overline{\Omega}}$. Thus , according to the resolvent equation, 
there is a point $a \in U$ belonging 
to the unbounded component of ${\bf C} \setminus {\overline{\Omega}}$ with the property that
the operator $({T^\ast} -{\overline{a}}{)^{-1}}$ leaves the subspace $K$ invariant. 
Then the operator $$
({T^\ast} -{\overline{b}}{)^{-1}}=({T^\ast}-{\overline{a}}{){-1}}(1+
({\overline{a}}-{\overline{b}})({T^\ast}-{\overline{a}}{)^{-1}} {)^{-1}} $$
still leaves $K$ invariant whenever $|b-a| \leq \| (T-a{)^{-1}}{\|^{-1}}$. 
By following a path disjoint of $\overline{\Omega}$ which joins $a$ to infinity
we find in finitely many steps that the operator $({T^\ast}-{\overline{c}}{)^{-1}}$
leaves the subspace $K$ invariant for $c$ in the neighbourhood of infinity. 
In conclusion, by taking the Taylor series of the resolvent function at infinity we
obtain that $K$ is a ${T^\ast}$-closed invariant subspace in the Hilbert space $H$
where $T$ acts. 

Let us denote $S=({T^\ast}{\mid _ K}{)^\ast}$, regarded as
an operator from $K$ to $K$. Since $F$ coincides with ${E_{\Omega}}$ on $(U \setminus
{\overline{\Omega}}{)^2}$, we have:
$$F(z,{\overline{w}}) = 1-\langle \Resw , \Resz \rangle ,\hspace{.2in}(z,w \in U \setminus {\overline{\Omega}}),$$
or equivalently $$
F(z,{\overline{w}}) = 1- \langle ({S^\ast} -{\overline{w}}{)^{-1}}\xi ,
({S^\ast}-{\overline{z}}{)^{-1}}\xi \rangle,\hspace{.2in}(z,w \in U \setminus 
{\overline{\Omega}}).$$
The difference between the two equations lies in the fact that the second one 
makes sense on the whole set $U \times U$.

Indeed, by the resolvent equation we obtain, for all $w,z,u \in U \setminus
{\overline{\Omega}}, z \neq u$:$$
\langle ({S^\ast}-{\overline{w}}{)^{-1}}\xi, ({S^\ast}-{\overline{z}}{)^{-1}}
({S^\ast}-{\overline{u}}{)^{-1}} \xi \rangle =$$
$$ (z-u{)^{-1}} [ 
\langle ({S^\ast} -{\overline{w}}{)^{-1}}\xi , ({S^\ast}-{\overline{z}}{)^{-1}}\xi \rangle
-\langle ({S^\ast} -{\overline{w}}{)^{-1}}\xi, ({S^\ast}-{\overline{u}}{)^{-1}}\xi \rangle ]=$$
$$ \frac{F(z,{\overline{w}})-F(u,{\overline{w}})}{z-u}.$$

This shows that for any vector $x \in K$ of the form 
$x=({S^\ast}-\overline{u}{)^{-1}}\xi$, the $K$-valued function
$ ({S^\ast}-{\overline{z}}{)^{-1}}x$ extends analytically
from $U \setminus {\overline{\Omega}}$ to $U$. However, the operator $S-z$  
may not be invertible 
on $K$, for all values $z \in U$. 

Since $$({T^\ast} -{\overline{z}})({S^\ast}-{\overline{z}}{)^{-1}}\xi = \xi,
\hspace{.2in}(z \in U),$$
we find $$
\| ({S^\ast}-{\overline{z}}{)^{-1}}\xi \| \geq \| \Resz \|,\hspace{.2in}(z \in U).$$
Consequently $$
F(z,{\overline{z}}) = 1-\| ({S^\ast}-{\overline{z}}{)^{-1}}\xi {\|^2} \leq 1-
\| \Resz {\| ^2}=0, \hspace{.2in}(z \in U \cap \Omega).$$

On the other hand, $$
F(z,{\overline{z}})={E_{\Omega}}(z,{\overline{z}})>0, \hspace{.2in}
(z \in U \setminus
{\overline{\Omega}}).$$

In conclusion $F(z,{\overline{z}})=0$ for every point $(\partial \Omega) \cap U$,
which shows that the set $(\partial \Omega) \cap U$ is included in a real analytic
subset of $U$.

It remains to invoke the local structure of real analytic sets (cf. for instance
[M] \S III.5.C) and to note that $\partial \Omega$ is the boundary of
an open set satisfying $\Omega = int{\overline{\Omega}}$. Indeed, for a small ball
$V$ centered at $\lambda$, the real analytic set $\{ z \in V; F(z,{\overline{z}})=0 \}$
consists of finitely many analytic arcs passing through $\lambda$
, some possibly being singular at $\Lambda$. These arcs divide $V$ into finitely many chambers,
each being included either in $\Omega$ or in ${\bf C} \setminus \Omega$. Then
$(\partial \Omega) \cap V$ consists of exactly those analytic semi-arcs starting
from $\lambda$ and dividing interior from exterior chambers. Thus the
set $(\partial \Omega) \cap V$ consists of finitely many full analytic semi-arcs
starting from the point $\lambda$ and ending on the boundary of the ball $V$.\bigskip

With the notation used in the proof of Proposition 3.1 we isolate the 
following result. For details concerning local spectral theory and the
significance of such a result in the context of abstract spectral decomposition
theories, we refer to [DS].

\paragraph{Corollary 3.2.}{\it Assume that the kernel ${E_{\Omega}}$ extends analytically
to a connected neighbourhood $U$ of a point $\lambda \in \partial \Omega$ and
$U$ intersects the unbounded connected component of ${\bf C} \setminus {\overline{\Omega}}$. 
Then the localized resolvent $({S^{\ast}}-\overline{z}{)^{-1}}x$ exists as an
analytic function in $z \in U$, for all vectors $x=({S^\ast}-\overline{u}{)^{-1}}\xi,
|u| \gg 0$.
}\bigskip

Later in this section we will see that in fact the unbounded connected component of
the complement of $\overline{\Omega}$ plays no special role. Any other component can
substitute it.

The following result is a partial converse to Proposition 3.1.

\paragraph{Theorem 3.3.} {\it Let $\Omega$ be a bounded planar domain and let
$\lambda \in \partial \Omega$ a point in whose neighbourhood $\partial \Omega$
is real analytic and smooth.

Then there is an open neighbourhood $U$ of $\lambda$ in ${\bf C}$, such that
the kernel ${E_{\Omega}}(z,{\overline{w}})$ extends analytically
from $z,w \in U \setminus {\overline{\Omega}}$ to $z,w \in U$.

Moreover, there is an invertible analytic function $f(z,{\overline{w}})$ 
defined for \\
$z,w \in U$, with the property that the analytic extension ${{E'}_{\Omega}}$ 
of the kernel ${E_{\Omega}}$ satisfies:
\begin{equation}
{{{E'}}_{\Omega}}(z,{\overline{w}})=(z-{w^\ast})f(z,{\overline{w}}), \hspace{.2in}
(z,w \in U),
\end{equation}
where $w \mapsto {w^\ast}$ is the Schwarz reflection function in $U \cap \partial \Omega$.}\bigskip

In the above statement it is implicit that we shrink $U$ to a neighbourhood of $\lambda$
where the anti-analytic Schwarz reflection map is defined, see [D],[Sh].\bigskip

{\it Proof.} We begin by recalling a simple computation from [P1], Scetion 5. Namely for the open 
unit disk $D \in {\bf C}$ ,we have: $$
{E_{D}}(z,{\overline{w}})= 1-\frac{1}{z{\overline{w}}}, \hspace{.2in}(z,w \in {\bf C} \setminus {\overline{D}}).$$

The idea of the proof is to transfer, via a conformal map, this identity 
to a neighbourhood in $\Omega$ 
of the smooth point $\lambda \in \partial \Omega$. The necessary computations are rather
long, but elementary. We also mention that in the following proof we do not use
the hyponormal operator attached to the domain $\Omega$.

Let $U$ denote an open ball centered at the point $\lambda$, such that $U \cap
\Omega$ is connected and simply connected and $(\partial \Omega) \cap U$ is a smooth
real analytic curve. Let $g:D \longrightarrow U \cap \Omega$ be a conformal map.
By Carath\'{e}odory Theorem and Schwarz Reflection Principle we can assume,
after possibly shrinking the radius of the ball U, that $g(1)=\lambda$, $g$ can be
extended conformally to a neighbourhood of the point $1 \in \partial D$ and hence
${g^{-1}}((\partial \Omega) \cap U) \subset \partial D$. Let $B$ be an open ball centered at
the point $1$, such that the conformal map $g$ is defined on $B$ and $U \subset g(B)$.

Since $$
{E_{\Omega}}(z,{\overline{w}})={E_{\Omega \cap U}}(z,{\overline{w}})
{E_{\Omega \setminus {\overline{U}} }}(z,\overline{w})\hspace{.2in}(z,w \in {\bf C}\setminus {\overline{\Omega}}),$$
and ${E_{\Omega \cap {\overline{U}} }}(z,{\overline{w}})$ is an analytic invertible function
on $U \times U$, it suffices to prove the statement for $\Omega \cap U$.

We change the variables as follows:$ \zeta=g(u), z=g(s), w=g(t)$. For $u \in B$ and
$s,t \in B \setminus \overline{D}$ we have:$$
{E_{\Omega \cap U}}(z,{\overline{w}}) = exp(\frac{-1}{2 \pi i}
\int_{\Omega \cap U} \frac{d{\overline{\zeta}}\wedge d\zeta}
                          {(\zeta -z)({\overline{\zeta}}-\overline{w})} )=$$
                          
$$ exp(\frac{-1}{2 \pi i} \int_{B \cap D} \frac
        { |g'(u){|^2}d{\overline{u}} \wedge du}
        {(g(u)-g(s))(\overline{g(u)}-\overline{g(t)})} )=   $$
  $$ exp(\frac{-1}{2 \pi i} \int_{B \cap U} h(u,{\overline{u}};s,{\overline{t}})
        \frac{d \overline{u} \wedge du}
        {(u-s)(\overline{u} -\overline{t})} ),$$
 where $$ 
 h(u,\overline{u}; s,\overline{t}) =\frac
 {|g'(u){|^2}(u-s)(\overline{u} -\overline{t})}
 {(g(u)-g(s))(\overline{g(u)}-\overline{g(t)})}.$$

 The function $h$ is analytic in all four variables, and by possibly shrinking
 the radius of the ball $B$ we can write:$$
 h(u,\overline{u};s,\overline{t})=1+(u-s){h_1}(u,s)+(\overline{u}-\overline{t})
 {h_2}(\overline{u},\overline{t})+(u-s)(\overline{u}-\overline{t}){h_3}(u,\overline{u};
 s,\overline{t}),$$
 where the new functions ${h_1},{h_2},{h_3}$ depend analytically on the
 respective variables (running through B).

 In what follows we replace $h$ by this additive decompositon in the above expression of
 ${E_{\Omega \cap U}}$. The first term $1$ in $h$ produces the factorization:
$$ {E_{B \cap D}}(s, \overline{t})=(1-\frac{1}{s\overline{t}}) {E_{D \setminus
 \overline{B} }}(s,\overline{t})\hspace{.2in}
 (s,t \in B \setminus \overline{D} ).$$

 For the second term, consider an analytic function ${H_1}(u,s)$ with the property that
$ \partial /\partial u {H_1}(u,s) ={h_1}(u,s)$ for $u,s \in B$. Then the corresponding
 integral (under the exponential function in the expression of ${E_{\Omega \cap U}}$) is:
$$ {I_1}(s,\overline{t}) = \frac{-1}{2 \pi i} \int_{B \cap D} 
\frac{{h_1}(u,s)}{\overline{u}-\overline{t}} d\overline{u} \wedge du=
\frac{1}{2 \pi i} \int_{\partial (B \cap D)} \frac{{H_1}(u,s)}{\overline{u}-
\overline{t}} d\overline{u} =$$
$$\frac{1}{2 \pi i} [\int_{D \cap \partial B} \frac{{H_1}(u,s)d\overline{u}}
{\overline{u} -\overline{t}} - \int_{B \cap \partial D} \frac{ {H_1}(u,s)du}
{u -{u^2}\overline{t}} ]=$$
$$\frac{1}{2 \pi i}\int_{D \cap \partial B} [
\frac{{H_1}(u,s)d\overline{u}}{\overline{u}-\overline{t}} -
\frac{ {H_1}(u,s)du}{u\overline{t}(u-{\overline{t} ^{-1}})} ] +
\frac{1}{2 \pi i} \int_{\partial(D \cap B)} \frac{ {H_1}(u,s)du}
{u\overline{t}(u-{\overline{t} ^{-1}})}.$$

In the last sum, the first intergral is analytic in $s,t$ belonging to a neighbourhood
of $1$ (such that ${\overline{t} ^{-1}} \in D \setminus \partial B$). The second integral
can be evaluated by Cauchy's Formula, for $t \in B \setminus \overline{D}$ and it is
equal to ${H_1}({\overline{t} ^{-1}},s)$. Thus, after shrinking the radius of the ball
B, if necessary, the function ${I_1}(s,\overline{t})$ extends
analytically from $B \setminus \overline{D}$ to $B$. 

Similarly one proves that the term ${h_2}$ produces an extendable analytic function,
from $s,t \in B \setminus \overline{D}$ to $s,t \in B$. The integral corresponding
to the term ${h_3}$ is obviously analytic in $B$. 

Summing up these computations, we have obtained an analytic function $A(s,\overline{t})$,
defined for $s,t \in B$, with the property that:$$
{E_{\Omega \cap U}} (g(s),\overline{g(t)}) = (1-\frac{1}{s\overline{t}}) exp A(s,\overline{t}),\hspace{.2in}(s,t \in B \setminus
\overline{D} ).$$
But the function $log(s)$ is well defined for $s \in B$, whence we can write:$$
{E_{\Omega \cap U}}(g(s),\overline{g(t)})= (s-{\overline{t} ^{-1}}) exp(A(s,\overline{t})-log(s)),$$
or changing back the variables:$$
{E_{\Omega \cap U}}(z,\overline{w})= ({g^{-1}}(z)-{g^{-1}}(w{)^\ast}) exp {A_1}(z,\overline{w}),
\hspace{.2in}(z,w \in U \setminus \overline{\Omega}),$$
where ${A_1}(z,\overline{w})$ is analytic in $z,w \in U$.

It remains to remark that:$$
{g^{-1}}(z)-{g^{-1}}(w{)^\ast}={g^{-1}}(z)-{g^{-1}}({w^\ast})=(z-{w^\ast})k(z,\overline{w}),$$
where the function $k(z,\overline{w})$ is invertible (and actually admits an analytic
logarithm) in the domain $z,w \in U$.

This completes the proof of Theorem 3.3.\bigskip

An analysis of the preceding proof shows that all computations make sense
before taking the exponential. In particular, let us remark that the function 
$log(z-{w^\ast})$ is well defined for $z,w \in U \setminus \overline{\Omega}$.
Thus we can state the following result.

\paragraph{Corollary 3.4.}{\it Let $\lambda \in \partial \Omega$ be a smooth,
real analytic point of the boundary of the domain $\Omega$. There exists an
open neighbourhood $U$ of $\lambda$ in ${\bf C}$ and an analytic function
$a(z,\overline{w})$ defined for $z,w \in U$, such that:\begin{equation}
\frac{-1}{2 \pi i} \int_{\Omega} \frac{d{\overline{\zeta}} \wedge d \zeta}
{(\zeta -z)(\overline{\zeta} -\overline{w}) } = log(z-{w^\ast})+a(z,\overline{w}),
\hspace{.2in}(z,w \in U \setminus \overline{\Omega}).
\end{equation} }\bigskip

In the opposite direction to Theorem 3.3 we give below an example of real
analytic singular boundary point across which the exponential kernel does not
extend analytically.

\paragraph{Proposition 3.5.} {\it Let $\Omega$ be a bounded domain of the complex
plane and let $\lambda \in \partial \Omega$ be a point in whose neighbourhood
the boundary $\partial \Omega$ is real analytic.

If an irreducible component of $\partial \Omega$ at $\lambda$ intersects
${\bf C} \setminus \overline{\Omega}$, then the kernel ${E_{\Omega}}(z,\overline{w})$
does not extend analytically to a full neighbourhood of the point 
$(z,w)=(\lambda ,\lambda)$. }\bigskip

{\it Proof.} Let $X$ be the germ at $\lambda$ of an analytic set, whose part of
half-branches define the germ of $\partial \Omega$ at $\lambda$. Let ${X_0}$ denote
an irreducible component of $X$ which intersects the set ${\bf C} \setminus \overline{\Omega}$.

Suppose that $F(z,\overline{w})$ is an analytic extension of ${E_{\Omega}}(z,\overline{w})$, from
a component of ${\bf C} \setminus \overline{\Omega}$ adjiacent to ${X_0}$
to a full neighbourhood of $(\lambda, \lambda)$ in ${{\bf C}^2}$. According to
Proposition 3.1, $F(z,\overline{z})=0$ for $z \in {X_0} \cap \partial \Omega$, hence for
$z$ belonging to a half-branch of ${X_0}$. 

Let $Y$ denote the germ of analytic set at $\lambda$ defined by the equation
$Y=\{ z; F(z, \overline{z})=0 \}.$ By elementary dimension theory (see for instance
[M] \S III.5), ${X_0} \cap Y$ is either $\{ \lambda \}$ or ${X_0}$. Since the first
alternative is excluded, we find that ${X_0} \subset Y$. But for a point 
$\mu \in {X_0} \setminus \overline{\Omega}$ we have:$$
0 < {E_{\Omega}}(z,\overline{z})=F(z,\overline{z})=0,$$
a contradiction!

Therefore the analytic set ${X_0}$ cannot intersect ${\bf C} \setminus
\overline{\Omega}$.\bigskip

We turn now to global analytic extensions of the exponential kernel of a domain. 
We prove below the extendability of the exponential kernel
along any continuous path on which the Schwarz function of a portion of 
the boundary extends.

Let $A$ be a smooth real analytic arc in the exterior boundary of the bounded
domain $\Omega$. Let $C$ denote the unbounded connected component of
${\bf C} \setminus \overline{\Omega}$.
Let $u(z)$ denote the Schwarz function of the arc $A$; $u$ is
an anaytic function in a neighbourhood of $A$ which coincides with $\overline{z}$
on $A$, see [D],[Sh]. 

Fix a point $\lambda \in A$ and consider a differentiable path 
without \\
self-intersection points, which joins $\lambda$ , within the domain 
$\Omega$, to a point 
$a \in \Omega$. Let $U$ be a tubular neighbourhood of $P$ in $\Omega \cup A$, such
that $U$ is a connected, simply connected domain with piecewise boundary, and
$A'=\partial U \cap \partial \Omega$ is an open analytic arc in $A$ , containing
of course the point $\lambda$.  We assume that the Schwarz function $u$ extends
analytically to $\overline{U}$. Thus the domain $C \cup A' \cup U$ is connected and simply
connected. (Visually, we attach the appendix $U$ to the simply connected component $C$).

Let $z,w \in C$ be fixed points. Then:
$$ \int_{\Omega} \inte = \int_{\Omega \setminus U} \inte + \int_{U} \inte .$$
Moreover, for a fixed determination of the logarithm ,defined ,say for 
$ \Re (w)>0,|w| > sup_{\zeta \in U} |u(\zeta)|,$ 
we obtain:
$$ \int_{U} \inte = \int_{\partial U} \binte =$$
$$\int_{A'} \binte +\int_{(\partial U) \setminus A'} \binte =$$
$$\int_{A'} \frac{log(u(\zeta)-\overline{w})d \zeta}
                 {\zeta -z} +\int_{(\partial U) \setminus A'} \binte =$$
                 $$ \int_{(\partial U)\setminus A'}[\binte-
\frac{log(u(\zeta)-\overline{w})d\zeta}{\zeta -w}].$$

These computations show that, for every $|w| \gg 0, \Re (w)>0,$ the function 
$z \mapsto {E_\Omega}(z,\overline{w})$ extends analytically from
$z \in C$ to $z \in C \cup A' \cup U$. Let us denote this extension by
$F(z,\overline{w})$.

At this moment the operator theory interpreation of the kernel ${E_{\Omega}}$
is again invoked to help. Let $T$ be the hyponormal operator attached to the domain $\Omega$.
In particular we know that:$$
F(z,\overline{w})= 1-\langle \Resw , \Resz \rangle ,\hspace{.2in}(z,w \in C).$$
Now we repeat the trick contained in the proof of Proposition 3.1. Let $$
K = \bigvee_{|w| \gg 0, \Re (w) >0} \Resw $$ be the ${T^\ast}$-invariant 
subspace generated by the vector $\xi$. Let ${S^\ast}$ be the restriction
of the operator ${T^\ast}$ to $K$. Then formula $$
F(z,w)=1-\langle ({S^\ast}-\overline{w}{)^{-1}}\xi,
({S^{\ast}}-\overline{z}{)^{-1}}\xi \rangle $$
holds for $|w| \gg 0, \Re (w) >0$ and $z \in C$.

Because the function $F(.,\overline{w})$ extends analytically to $C \cup A'
\cup U$, we infer, exactly as in the proof of Proposition 3.1, that the function
$F(z,\overline{w})$ extends analytically from $z,w \in C$ to $z,w \in C \cup A' 
\cup U$. 

Soon we will consider only domains with real algebraic boundary, so that the
above suppositions are not excessive.

In conclusion we have proved the following result.

\paragraph{Theorem 3.7.} {\it Let $\Omega$ be a bounded domain and let
$A$ be an analytic arc in the exterior boundary of $\Omega$. If the Schwarz
function of the arc $A$ extends analytically in the interior of $\Omega$ to an open set
$U \subset \Omega$, then the kernel ${E_{\Omega}}(z,\overline{w})$ extends
analytically from the unbounded component $C$ of ${\bf C} \setminus \overline{\Omega}$
to $z,w \in C \cup A \cup U$.}\bigskip

There are examples which show that, to different analytic continuations of the 
Schwarz function, correspond in general different analytic continuations of the
exponential kernel.

In order to extend the results of this section to other connected components of the
comlement of $\Omega$, it is sufficient to see how the kernel ${E_{\Omega}}$
changes under linear fractional transformations.

\paragraph{Proposition 3.8.}{\it Let $\Omega$ be a bounded planar domain
and let $a,b$ be complex numbers, $b \neq 0$. Suppose that $0 \in {\bf C}
\setminus \overline{\Omega}$. Then for every pair $z,w \in {\bf C}\setminus 
\overline{\Omega}$:\bigskip

a). ${E_{\Omega^{-1}}}({z^{-1}},{{\overline{w}}^{-1}})=
\frac{ {E_\Omega}(0,0) {E_\Omega}(z,\overline{w})}
     { {E_\Omega}(0,\overline{w}) {E_\Omega}(z,0) }$;\bigskip

b). ${E_{\Omega +a}}(z+a,\overline{w+a})= {E_\Omega}(z,\overline{w})$;\bigskip

c). ${E_{ b \Omega }}(bz, \overline{bw})= {E_\Omega}(z,\overline{w})$.}\bigskip

We have denoted ${\Omega^{-1}}=\{ {z^{-1}}; z \in \Omega \}$ and so on.
The proof of Proposition 3.8 consists in simple changes of variables in
the integral defining the exponential kernel. The details are left to the reader.

It is important to remark that the denominator in formula $a)$ above does not vanish.
Indeed, $$
{E_{\Omega}}(z,0) = 1 -\langle {T^{\ast -1}}\xi , \Resz \rangle, $$
and $max\{ \| {T^{\ast -1}}\xi \|, \| \Resz \| \} <1$ for $z \in {\bf C} \setminus
\overline{\Omega}$.

\section{Domains with real algebraic boundary}

In view of the general convexity results of Section 2, the extremal solutions
of the $L$-problem of moments are among the characteristic functions of 
semialgebraic domains. In this section we begin a study of some specific 
properties of these domains.

We start with a couple of examples. Let $0<r <R$ be fixed real numbers and let
$A(r,R)$ be the annulus $A(r,R)=\{ z \in {\bf C}; r<|z|<R \}.$  Proposition 3.8
shows that: 
\begin{equation}
{E_{A(r,R)}}(z,\overline{w})= \left\{ \begin{array}{cc}
   \frac{1-{R^2}z \overline{w}}{1-{r^2}z \overline{w}} & |z|,|w| >R \\
   \frac{1-{r^2}z \overline{w}}{1-{R^2}z \overline{w}} & |z|,|w| <r.
   \end{array}
   \right.
   \end{equation}

In the second example we do not compute explicitly the exponential kernel,
but instead we locate the singularities of the analytic extensions of it.
Take for instance a non-degenerated triangle $T$ in
$ {\bf C}$. Then, according to Theorem 3.7, the kernel ${E_T}$
extends analytically from the complement of $T$ to the complement of the
union of two arbitraily chosen sides of $T$. Indeed, in this case the Schwarz function
of each side is entire, see [D] Chapter 10. As a consequence, formula $$
\int_{T} f dA =u(f),\hspace{.2in}(f \in {\cal O}(\overline{T})),$$
holds, where $u$ is an analytic functional (that is a linear continuous functional
on  ${\cal O}(\overline{T})$) carried by the respective two sides of $T$. 
(The reader can actually find explicitly the functional $u$ in the
above quadrature formula).For the distinction between the carrier of an analytic functional and the
support of a distribution we refer to [H] Section IX.9.1 or [Ma].

Next we prove that similar phenomena occur for arbitrary semi-algebraic domains.
Let $\Omega$ be a bounded planar domain with real algebraic boundary. Let $C$
denote the unbounded component of ${\bf C} \setminus \overline{\Omega}$ and let
${C_1},\ldots ,{C_n}$ be the bounded connected components of the same set.

The 
Schwarz function of each analytic arc contained in the boundary of $\Omega$
is an algebraic function. Let $A$ denote the set of the ramification points 
in $\overline{\Omega}$ of all these local Schwarz functions, union with the
singular or multiple points of the boundary of $\Omega$. Thus, $A$ is a finite
subset of $\overline{\Omega}$. We claim that there is a finite set $E$ of
continuous curves contained in $\overline{\Omega}$ and passing through all
points of the set $A$, such that $\Omega '=(\overline{C} \cup \Omega) \setminus
E$ is a connected open set and the function $\overline{z}$ extends analytically
from $\partial \Omega \setminus E$ to $\Omega \setminus E$. The reader can 
either consult [AG] Chapter V for the existence
of such continuous cuttings or prove directly the claim.

With this preparation and Theorem 3.7 we can state the following result.

\paragraph{Theorem 4.1.}{\it Let $\Omega$ be a bounded planar domain with
real-algebraic boundary
and let $\Omega '$ be the domain consisting of the unbounded component of the
complement of $\overline{\Omega}$ union with $\Omega$, minus the cuttings explaind above.

Then the kernel ${E_{\Omega}}(z,\overline{w})$ extends analytically from
$z,w$ belonging to the unbounded component of the complement of $\overline{\Omega}$ to $
z,w \in \Omega '$.}\bigskip

As the preceding examples show, the kernel ${E_{\Omega}}$ may extend analytically across
some of the curves in the set $E$, but in general these extensions cannot be 
glued.

\paragraph{Corollary 4.2.} {\it With the notation in Theorem 4.1 we
obtain the following quadrature identity: $$
\int_{\Omega} f dA = u(f), \hspace{.2in}(f \in {\cal O}(\overline{\Omega})),$$
where $u$ is an analytic functional carried by $\partial \Omega '$.}\bigskip

{\it Proof.} Indeed, let $u(z)$ be the analytic function which extends 
$\overline{z}$ from $\partial C \cap \Omega '$ to $\Omega '$.
Then for an analytic function in a neighbourhood of $\overline{\Omega}$ we
obtain :$$
\frac{1}{2i} \int_{\Omega} fdA = \int_{\partial \Omega} f(z)\overline{z} dz=$$
$$
\sum_{j=1}^{n} \int_{C_j} f(z)\overline{z} dz +\int_{\partial C \cap E} f(z)
\overline{z} dz + \int_{\partial C \setminus E} f(z) u(z) dz.$$

By Cauchy's Theorem we can replace the integration curve in the last integral
by a curve which is arbitrarily close to $E$ union with the interior boundaries 
of $\Omega$ and Corollary 4.2 follows. \bigskip

According to Proposition 3.8 similar results hold for any connected component of the
complement of the domain $\Omega$. However, the first example above shows that
the analytic extensions of the exponential kernel across different connecetd components
of the boundary of $\Omega$ do not coincide.

As a final remark we return to the original unsolved question which 
has motivated the whole paper . Namely, from Section 2
we know that the exponential kernels of the extremal domains exhibited there 
are finitely determined by their Taylor coefficients
at infinity, up to a certain prescribed degree; on the other hand Theorem 4.1 above shows that
they have an analytic continuation configuration consisting essentially 
of finitely many ramification points. The question is whether these data 
suffice for finding
some closed expression for the exponential kernels of these extremal domains 
(involving most likely some special functions).\bigskip

\section{ A characterization of the exponential kernel}

To reveal another face of the same exponential kernel, we derive below
a characterization of it in terms of positive definite functions. The principle, 
of obtaining such a characterization in terms of postive definite kernels is not new. 
It goes back to the work of de Branges on Hilbert spaces of analytic 
functions [dB]; the same technique was exploited later by Pincus and Rovnyak
[PR], Carey and Pincus [CP] and several other authors whose interests came in 
contact with determining or charactersitic functions of various classes of operators.
For another notable example see also Liv\v{s}ic [L].
In spite of several close similarities with the mentioned works, we believe
that the details contained in this section are new.

To simplify notation we put $\hat{\bf C}={\bf C} \cup \{ \infty \},
{\bf D}=\{ z \in {\bf C}; |z| <1 \}$ and ${{\bf C}^\ast}={\bf C}\setminus
\{ 0 \}$. For a measurable function $g:{\bf D} \longrightarrow [0,1]$  we 
denote: $$
{E_g}(z,\overline{w})= exp (\frac{-1}{\pi}\int_{\bf D}
\frac{g(\zeta)}{ (\zeta -z)(\overline{\zeta} -\overline{w})} dA(\zeta)) ,
\hspace{.2in}(|z|,|w|>1).$$ Let us also remark that the above function
extends analytically to a function $$
{E_g}: (\hat{\bf C} \setminus \overline{\bf D}{)^2} \longrightarrow {{\bf C}^\ast}.$$
       
The question we address in this section is to find a set of minimal conditions
which characterizes an analytic function \\
$E:({\bf C}\setminus \overline{\bf D}
{)^2} \longrightarrow {{\bf C}^\ast}$ to be of the form $E={E_g}$ for a 
measurable function $g$ as above.

One obvious condition is:\begin{equation}
E(\infty ,\overline{w})=E(z,\infty)=1 ,\hspace{.2in}(z,w \in \hat{\bf C}
\setminus \overline{\bf D}).\end{equation}

In order to state the next condition we define a new kernel \\
$F:(\hat{\bf C} \setminus \overline{\bf D}{)^
4}\longrightarrow {\bf C}$
by the formula:\begin{equation}
F({z_1},\overline{z_2};{w_1},\overline{w_2})=
\frac{ E({z_1},\overline{w_2})E({w_1}, \overline{z_2})-E({z_1},\overline{z_2})
E({w_1},\overline{w_2})}
{({w_1}-{z_1})(\overline{w_2}-\overline{z_2})E({z_1},\overline{w_2})}.
\end{equation}

Whenever we encounter an analytic function $h(z)$, the quotient $$
\frac{h(z)-h(w)}{z-w},\hspace{.2in}(z \neq w),$$
is extended analytically across the diagonal $(z=w)$ by the value $h'(z)$.As a
matter of terminology the inequality $K({z_1},\overline{z_2};{w_1},
\overline{w_2}) \succ 0$ means that the kernel $K$ is non-negatively definite,
that is $$
\sum_{k,l=1}^{N}  K({s_k},\overline{t_k};{t_l},\overline{s_l}){\lambda _k}\overline
{\lambda _l} \geq 0,$$
for every finite set of points $\{ ({s_k},{t_k}), 1 \leq k \leq N \}$ in the
domain of $K$ and
every complex numbers ${\lambda _k}, 1 \leq k \leq N.$\bigskip
 
\paragraph{Theorem 5.1.} {\it Let $E:(\hat{\bf C} \setminus \overline{\bf D}{)^2
}\longrightarrow {{\bf C}^\ast}$ be an analytic function with the normalization
property (9) and let $F$ be the associated kernel (10).

There is a measurable function $g:{\bf D} \longrightarrow [0,1]$ with the
property that $E={E_g}$ if and only if $$
F({z_1},\overline{z_2};{w_1},\overline{w_2}) \succ
{z_1}\overline{w_2}F(\variabile)-({z_1}\overline{w_2}F(\variabile){)_{{z_1}=\infty}} -$$ \begin{equation}
({z_1}\overline{w_2}F(\variabile){)_{{w_2}=\infty}} +({z_1}\overline{w_2}
F(\variabile){)_{ {z_1}={w_2}=\infty}} \succ 0.
\end{equation} }\bigskip

Note that the second term in the above positivity condition is a second order
difference at infinity of the function $F$.\bigskip

{\it Proof. The necessity.} For this part of the proof we use the known 
factorization of the kernel $1-{E_g}$ in terms of the associated hyponormal
operator with rank one self-commutator, see [P1],[P2].

Let $g$ be a measurable function as in the statement and let $T$ denote the
irreducible hyponormal operator with rank one self-commutator which has the
principal function equal to $g$, almost everywhere. Let us denote (as before)
$[ {T^\ast},T]=\xi \otimes \xi$ and let us recall the basic formula:
\begin{equation}
{E_g}(z,\overline{w})=1-\langle \Resw ,\Resz \rangle .
\end{equation}
Note that we have tacitly made the normalization $supp(g) \subset \overline
{\bf D}$, hence $\sigma (T) \subset \overline{\bf D}$ and $\| T \| \leq 1$,
(because $T$ is a hyponormal operator, see [MP], Corollary 3.1.4).

Next we need a few elementary identities with resolvents,all stated for the
current variables $u,v,{z_1},\ldots$ outside the closed unit disk:
\begin{equation}
\presu \presv =\frac{\presu -\presv}{u-v}
\end{equation}
and
\begin{equation}
\resu \presv = \presv \resu + 
\end{equation}
$$\resu \presv (\xi \otimes \xi) \presv \resu.
$$

In particular, from formula (14) we obtain:$$
\langle \resu \presv \xi, \xi \rangle =\langle \presv \resu \xi, \xi \rangle +$$
$$\langle \resu \presv \xi, \xi \rangle \angle \presv \resu \xi, \xi \rangle,$$
or equivalently:$$
(1-\langle \presv \resu \xi, \xi \rangle)(1+\langle \resu \presv \xi, \xi \rangle)=1,$$
that is:
\begin{equation}
1+\langle \resu \presv \xi ,\xi \rangle =\frac{1}{{E_g}(v,\overline{u})}.
\end{equation}

We claim that:\begin{equation}
F(\variabile)=\langle \presza \reszb \xi, \preswb \reswa \xi \rangle.
\end{equation}

Indeed, according to these identities we obtain (denoting $E={E_g})$:$$
\langle \presza \reszb \xi, \preswb \reswa \xi \rangle =$$
$$ \langle \PRESwa \RESwb \presza \reszb \xi, \xi \rangle =$$
$$ \langle \PRESwa \presza \RESwb \reszb \xi, \xi \rangle +$$
$$\langle \PRESwa \RESwb \presza \xi, \xi \rangle \langle \presza \RESwb
\reszb \xi, \xi \rangle =$$
$$ \frac{\langle \PRESwa \RESwb \reszb \xi, \xi \rangle -
\langle \presza \reszb \xi, \xi \rangle}{\overline{w_2}-\overline{z_2}}=$$
$$\frac{\langle\PRESwa \RESwb \xi, \xi \rangle -
\langle \PRESwa \reszb \xi, \xi \rangle}{({w_1}-{z_1})(\overline{w_2}-
\overline{z_2})}+$$
$$\frac{\langle \presza \reszb \xi, \xi \rangle -
\langle \presza \RESwb \xi, \xi \rangle}
{({w_1}-{z_1})(\overline{w_2}-\overline{z_2})}+$$
$$( \langle \PRESwa \presza \RESwb \xi, \xi \rangle +$$
$$\langle \PRESwa \RESwb \presza \xi, \xi \rangle
\langle \presza \RESwb \xi, \xi \rangle)\times $$
$$\frac{\langle \presza \RESwb \xi, \xi \rangle -
\langle \presza \reszb \xi, \xi \rangle}
{\overline{w_2}-\overline{z_2}}=$$
$$\frac{E({w_1},\overline{z_2})+E({z_1},\overline{w_2})-E({z_1},\overline{z_2})
-E({w_1},\overline{w_2})}
{({w_1}-{z_1})(\overline{w_2}-\overline{z_2})}+$$
$$\frac{\langle \PRESwa \presza \RESwb \xi, \xi \rangle(E({z_1},\overline{z_2})-
E({z_1},\overline{w_2})}
{E({z_1},\overline{w_2})(\overline{w_2}-\overline{z_2})}=$$
$$\frac{E({w_1},\overline{z_2})+E({z_1},\overline{w_2})-E({z_1},\overline{z_2})
-E({w_1},\overline{w_2})}
{({w_1}-{z_1})(\overline{w_2}-\overline{z_2})}+$$
$$\frac{(E({z_1},\overline{w_2})-E({w_1},\overline{w_2}))
(E({z_1},\overline{z_2})-E({z_1},\overline{w_2}))}
{E({z_1},\overline{w_2})({w_1}-{z_1})(\overline{w_2}-\overline{z_2})}=$$
$$ \frac{ E({z_1},\overline{w_2})E({w_1},\overline{z_2})-
E({w_1},\overline{w_2})E({z_1},\overline{z_2})}
{E({z_1},\overline{w_2})({w_1}-{z_1})(\overline{w_2}-\overline{z_2})}.$$

Thus relation (16) is verified. It remains to remark that :$$
T \presza \reszb \xi = $$
$$ \reszb \xi +{z_1}\presza \reszb \xi =$$
$$ {z_1} \presza \reszb \xi -({z_1}\presza \reszb \xi {)_{ {z_1}=\infty}}.$$

In conclusion the
 positivity conditions (11) in the statement become:
$$0 \prec \langle T \presza \reszb \xi, T \preswb \reswa \xi \rangle \prec $$
$$\prec \langle \presza \reszb \xi, \preswb \reswa \xi \rangle. $$

Since $T$ is a contraction these two positivity conditions are evidently
true.\bigskip

{\it The sufficiency.} Let $E$ be an analytic function which satisfies the 
normalization and positivity conditions in the statement. We want to prove that
$E={E_g}$, where $g:{\bf D} \longrightarrow [0,1]$ is a measurable function.
This in turn is equivalent in finding a hyponormal operator $T$ with rank-one
self-commutator $[{T^\ast},T]=\xi \otimes \xi$ which has $g$ as pricipal 
function and hence $E$ as associated determinantal function.

Since the kernel $F$ was supposed to be non-negatively definite, Kolmogorov's
factorization theorem implies the existence of a separable, complex
Hilbert space $H$ and an $H$-valued analytic function:$$
\rho : (\hat{\bf C} \setminus \overline{\bf D}{)^2}\longrightarrow H,$$
such that:\begin{equation}
F(\variabile )=\langle \rho ({z_1},\overline{z_2}), \rho({w_2},\overline{w_1})
\rangle,\hspace{.2in}({z_j},{w_j} \in \hat{\bf C} \setminus \overline{\bf D};
j=1,2).
\end{equation}
In addition, we can assume without loss of generality that the image of the
function $\rho$ spans the whole Hilbert space $H$.

By assumption, $F(\infty, \overline{z_2}; {w_1}, \overline{w_2})=
F({z_1},\infty ; {w_1},\overline{w_2})=0$, therefore:$$
\rho(\infty, \overline{z})=\rho(z,\infty)=0,\hspace{.2in}(z \in \hat{\bf C} \setminus
\overline{\bf D}).$$ In particular, both limits $lim_{{z_1} \rightarrow \infty}
{z_1}\rho({z_1},\overline{z_2})$ and $ lim_{{z_2} \rightarrow \infty}
\overline{z_2}\rho({z_1},\overline{z_2})$ exist.

We define a linear transformation on the range of $\rho$ by the formula:
\begin{equation}
T\rho({z_1},\overline{z_2})={z_1}\rho({z_1},\overline{z_2})-
({z_1}\rho({z_1},\overline{z_2}){)_{{z_1}=\infty}}.
\end{equation}

Let $n$ be a positive integer and let us choose arbitrary elements \\ 
${z_1}(k),{z_2}(k) \in \hat{\bf C}
\setminus \overline{\bf D},{\lambda _k}\in {\bf C}, 
1 \leq k \leq n.$
In view of condition (11) in the statement of Theorem A.1 we have:$$
\| T \sum_{k=1}^{n} {\lambda _k}\rho({z_1}(k),\overline{ {z_2}(k)}) \|
\leq \| \sum_{k=1}^{n} {\lambda _k}\rho({z_1}(k),\overline{ {z_2}(k)}) \|.$$
Therefore, the map $T$ extends linearly to a contraction defined on the whole
space $H$. We denote its extension by the same symbol $T$.

Our next aim is to prove the formula:\begin{equation}
{T^\ast}\rho({z_1},\overline{z_2})=\overline{z_2}\rho({z_1},\overline{z_2})-
E({z_1},\overline{z_2})(\overline{z_2}\rho({z_1},\overline{z_2}){)_{ {z_2}=
\infty}}.
\end{equation}

To this end we choose arbitrary points ${z_1},{z_2},{w_1},{w_2} \in
\hat{\bf C} \setminus \overline{\bf D}$ and compute:
 $$ \langle {T^\ast}\rho({z_1},\overline{z_2}),\rho({w_2},\overline{w_1})\rangle
-\langle\rho({z_1},\overline{z_2}),T\rho({w_2},\overline{w_1})\rangle=$$
$$\langle\overline{z_2} \rho({z_1},\overline{z_2}),\rho({w_2},\overline{w_1})
\rangle - E({z_1},\overline{z_2})(\langle \overline{z_2} \rho({z_1},
\overline{z_2}), \rho({w_2},\overline{w_1})\rangle{)_{ {z_2}=\infty}}-$$
$$\overline{w_2}\langle \rho({z_1},\overline{z_2}),\rho({w_2},\overline{w_1})
\rangle + (\langle \rho({z_1},\overline{z_2}),{w_2}\rho({w_2},\overline{w_1})
\rangle {)_{ {w_2}=\infty}}=$$
$$(\overline{z_2}-\overline{w_2})F(\variabile)-E({z_1},\overline{z_2})
(\overline{z_2}F(\variabile){)_{ {z_2}=\infty}}+$$
$$(\overline{w_2}F(\variabile)
{)_{ {w_2}=\infty}}=$$
$$  \frac{- E({z_1},\overline{w_2})E({w_1},\overline{z_2})+E({z_1},
\overline{z_2})E({w_1},\overline{w_2}) +E({z_1},\overline{z_2})(E({z_1},
\overline{w_2})-E({w_1},\overline{w_2}))}
{({w_1}-{z_1})E({z_1},\overline{w_2})}+$$
$$\frac{E({w_1},\overline{z_2})-E({z_1},\overline{z_2})}
{ {w_1}-{z_1}}=0.$$

Thus formula (19) is verified.

Let us remark that, denoting:$$
\xi = ({z_1}\overline{z_2}\rho({z_1},\overline{z_2}){)_{ {z_1}={z_2}=\infty}},
$$ we have:$$
({T^\ast}-\overline{z_2})(T-{z_1})\rho({z_1},\overline{z_2})=
({T^\ast}-\overline{z_2})(-{z_1}\rho({z_1},\overline{z_2}){)_{ {z_1}=\infty}}=$$
$$( E({z_1},\overline{z_2})({z_1}\overline{z_2}\rho({z_1},\overline{z_2})
{)_{ {z_2}=\infty}}{)_{ {z_1}=\infty}}=\xi,$$
(because $E(\infty, \overline{z_2})=1$). Whence we find the formula:
\begin{equation}
\rho({z_1},\overline{z_2})= \presza \reszb \xi, \hspace{.2in}({z_1},{z_2}\in
\hat{\bf C} \setminus \overline{\bf D}).
\end{equation}
Consequently we obtain:$$
[{T^\ast},T]\rho({z_1},\overline{z_2})=[{T^\ast}-\overline{z_2},
T-{z_1}]\presza \reszb \xi=$$
$$ \xi -E({z_1},\overline{z_2})\xi=(1-E({z_1},\overline{z_2}))\xi.$$
Therefore the operator $T$ has rank-one self-commutator, and the vector
$\xi$ spans the range of $[{T^\ast},T]$.

Finally we return to formula (10) and remark that:$$
1-E({z_1},\overline{z_2})=({w_1}\overline{w_2}F(\variabile ){)_
{{w_1}={w_2}=\infty}}=\langle \rho({z_1},\overline{z_2}),\xi \rangle,$$
so that $$
[{T^\ast},T]\rho({z_1},\overline{z_2})=\langle \rho({z_1},\overline{z_2})\xi,
\xi \rangle \xi.$$
This proves that $[{T^\ast},T]=\xi \otimes \xi $.

In conclusion $E={E_g}$, where $g$ is the principal function of the operator
$T$. 

This finishes the proof of Theorem 5.1.

\paragraph{Remark 5.2.} By changing the variables ${u_j}=\frac{1}{z_j},
{v_j}=\frac{1}{w_j}, j=1,2,$ we can define the function:$$
G({u_1},\overline{u_2}; {v_1},\overline{v_2})=
\frac{F(\variabile)}{{z_1}\overline{w_2}},$$
so that $G$ is analytic in the polydisk ${{\bf D}^4}$.

For an analytic function $h(z), z \in {\bf D},$ we define the difference
of $h$ at zero by:$$
{\Delta _z}h(z)=\frac{h(z)-h(0)}{z}.$$
Then condition (11) becomes:
\begin{equation}
0 \prec {\Delta _{u_1}}{\Delta _{ \overline{v_2}}}G({u_1},\overline{u_2};
{v_1},\overline{v_2}) \prec G({u_1},\overline{u_2};{v_1},\overline{v_2}).
\end{equation}.

Let ${\cal G}$ denote the class of all analytic functions $G :{{\bf D}^4}
\longrightarrow {\bf C}$ which satisfy the positivity conditions (21) and
have the structure derived from fomula (10), where $E$ is subject to the
normalization (9). 

Then the truncated, or full, $L$-problem of moments treated in
the previous sections (and in the papers [P1],[P2]) is equivalent to 
the following interpolation  problem
for the class ${\cal G}$:$$G \in {\cal G}$$ and $$
({{(\partial /\partial{u_1}})^m}{{(\partial /\partial{\overline{u_2}}})^n}
G)(0,0;0,0)={b_{mn}},\hspace{.2in}(0 \leq m \leq m+n \leq N).$$

This is a two-dimensional variant of the classical Carath\'{e}odory-Fej\'{e}r
problem (see for instance [FF]). On the basis of our previous results obtained for the $L$-problem of moments, 
we know for the the above interpolation problem how to describe its solvability 
in positivity terms (in the case $N=\infty$), while for the corresponding truncated interpolation problem we know a description
of all its extremal solutions (for $N$ finite). In view of the bijection between
the class of functions ${\cal G}$ and the measurable functions $g:{\bf D}
\longrightarrow [0,1]$, the class ${\cal G}$ has a natural convex structure
hidden in the free parameter $g$.

\end{document}